\documentclass[a4paper,10pt]{article}
\usepackage[utf8x]{inputenc}
\usepackage[russian]{babel}
\usepackage{amssymb}

\title{Special Bohr - Sommerfeld geometry on Riemann surfaces: toy problems.}
\author{Nikolay A. Tyurin\footnote{BLTPh JINR (Dubna) and  NRU HSE (Moscow), {\bf ntyurin@theor.jinr.ru}. {\it The author  was supported by RSF grant, project 14-21-00053 dated 11.08.14.}}}

\begin{document}

\maketitle

In paper [1] one proposes some new construction in lagrangian geometry of algebraic varieties. Every algebraic variety being endowed with a Kahler metric
of the Hodge type can be regarded as a symplectic manifold, and it is reasonable to study its lagrangian submanifolds. By the definition given in [1]
this new construction is mostly applicable to the case of simply connected manifolds, otherwise the notion of Bohr - Sommerfeld lagrangian submanifold
is not intrinsinc and for possible extensions one needs certain additional choice of the prequantization connection. But since we are speaking about
algebraic varieties it could happen that certain connection is already fixed in our picture so for this case quite natural to introduce the notion of Bohr - Sommerfeld
w.r.t. this particular connection and therefore to study special Bohr - Sommerfeld lagrangian submanifolds.

Briefly recall, let $(M, \omega)$ --- compact simply connected symplectic manifold s.t. the cohomology class of the symplectic form $\omega$ is integer. The prequantization data
$(L, a)$ consists of a line bundle and a hermitian connection on it such that the curvature form $F_a = 2 \pi  i \omega$. Then a lagrangian submanifold $S \subset M$ is said to be
Bohr - Sommerfeld (BS for short) iff the restriction $(L, a)|_S$ admits covariantly constant sections. If the corresponding  covariantly constant section is denoted as $\sigma_S$ then
$S$ is said to be special Bohr - Sommerfeld lagrangian w.r.t. a section $s \in \Gamma(M, L)$ iff the restriction $s|_S$ is presented as $s|_S = e^{i c} f \sigma_S$
where $c \in \mathbb{R}$ and $f$ is positive real function on $S$ (details see in [1]). The main case here is when $s$ is not just a smooth section but is holomorphic
section w.r.t. to certain compatible integrable complex structure $I$ on $M$. This means that we are speaking about an algebraic variety since $(M, \omega, I)$
is algebraic. 

If $M$ is not simply connected then the hermitian connection $a$ is not uniquelly defined up to gauge transformations therefore to define BS lagrangian submanifolds we must
fix certain $a4$, but it could be done automatically if $a$ is fixed from the data $(\omega, I)$: let the determinant Levi - Civita connection $a_{LC}$ on $det T^*M$
yields the prequantization condition $F_{a_{LC}} = 2 \pi k \omega$, where $k$ is a constant, then $a_{LC}$ can be taken for the definition of BS lagrangian submanifolds.
It's clear that the prequantization condition is exactly the Kahler - Einstien equation thus SBS geometry is still aplicable to the case of Fano varieties, Calaby - Yau
varieties etc.

Here we present the simplest case of ``etc'' varieties namely of Riemann surfaces, posting the simplest questions which arise in this way.

Let $\Sigma$ be a Riemann surface of genus $g> 1$ equipped with a fixed complex structure $I$. Then this complex structure can be extended
to the Kahler triple  $(G, I, \Omega)$, where $G$ is a riemannian metric of constant negative curvature compatible with $I$. Up to constant this riemannian metric is unique,
and we normalize it by the condition $\int_{\Sigma} \Omega = 2g-2$. It is well known that this metric is a solution of the Kahler - Einstein equation (see, f.e. [2]).

The Kahler structure $(G, I, \Omega)$ induces the corresponding hermitian structure on the complex line bundle $T^* \Sigma$, which we denote as $K_{\Sigma}$
and call the canonical bundle following the algebro - geometrical traditions. In the presence of the complex structure one has a finite dimensional
subspace $H^0(\Sigma, K_{\Sigma})$ of holomorphic section in the big space of all smooth sections; we call the sections holomorphic differentials. 

Any holomorphic differential  $\rho \in H^0(\Sigma, K_{\Sigma})$ has generically $2g-2$ zeros $P_1, ..., P_m$ but in particular cases certain $P_i$
are multiple zeros, in which case $m < 2g-2$. Consider the following smooth real function 
$$
\Psi(\rho) = - \rm{ln} \vert \rho \vert_h,
$$
which is correctly defined on the punctured surface $\Sigma \backslash \{P_1, ..., P_m \}$, where $P_i$ --- zeros of holomorphic differential $\rho$,
and the norm is taken with respect to the hermitian structure on the canonical bundle. Since the function $\Psi(\rho)$ is strictly convex w.r.t. the complex structure
$I$ (see Proposition 6 in [1]) therefore all its ``finite'' isolated  critical points have Morse indecies 0 or 1 (details see f.e. in [3]). Thus we have only minima and saddle
points. 

It's not hard to see, that for a generic holomorphic differential the corresponding function $\Psi (\rho)$ has non degenerated isolated critical points only,
and the number of these points is finite. This remark implies the first natural question: for generic holomorphic differential $\rho$ estimate 
(or find) the number of minima for the function $\Psi(\rho)$ on $\Sigma \backslash \{p_1, ..., p_m \}$.

Then the number  of saddle points is given from the formula
$$
\sharp {\rm min }  - \sharp {\rm saddle}  + 2g-2 = 2 - 2g
$$
for the generic case $m = 2g-2$.  Of course, the number varies along the projective space 
$\mathbb{C} \mathbb{P}^{g-1} = \mathbb{P} H^0(\Sigma, K_{\sigma})$, since for certain holomorphic differentials zeros are mutiple etc., the but what are the possible values?  

Furthemore, following the SBS strategy we are interested in the base set $B_{\rho} \subset \Sigma$, defined in [1], Definition 3. It can be defined in our present situation:
take all finite critical points and all finite trajectories of the gradient flow of $\Psi(\rho)$, this union is the base set $B_{\rho}$.

Note that despite of the fact that for certain $\rho$ the finite critical points of $\Psi(\rho)$ may form 1 - dimensional subsets in $\Sigma$,
the base set $B_{\rho}$ is still 1 - dimensional. Indeed, the critical subsets cann't be 2 - dimensional since $\Psi(\rho)$ is Kahler potential
for $\Omega$, on the other hand $B_{\rho}$ by the very definition must be stable w.r.t. to the gradient flow, and if it contains 2 - dimensional
components then this component must be transported by the gradient flow to a saddle point, but then this saddle point should have
2 - dimensional negative subspace in the tangent space which is impossible. 

Moreover, SBS geometry can help to prove here that the number of finite gradient trajectories for generic $\rho$ is {\bf finite}. Indeed, we just modify the arguments of Proposition 2 in [1]:
if two finite gradient trajectories with the same ends are homotopical one to each other then they must bound a domain whose symplectic area equals to an integer  multiple of $2g-2$,
but since the total sysmplectic area of $\Sigma$ is $2g-2$ it is impossible. On the other hand two gradient trajectories cann't intersect each other, consequently
we can get only finite number of finite trajectories.

Note that if we oversee the picture on $\Sigma$ we get that $B_{\rho}$ looks like a graph where all minima present the vertices and the saddle points are not visible, being
just marked points on the edges of the graph. Indeed, two trajectories meet at a saddle point $p$ such that they a tangent to the negative subspace $T^-_p \Sigma \subset T_p \Sigma$,
but since $T^-_p\Sigma$ is 1 - dimensional these two trajectories give us one smooth path with one marked point $p$. Going along an edge we can meet several
marked points since in general it is possible to have several saddle points on the same line combinig the finite trajectories. Moreover, each part
of each edge is naturally oriented since the function $\Psi(\rho)$ changes along these segments, and we can attach to each segment the approriate sign.  
Note that every edge must have at least one marked point.

Therefore for a generic holomorphic differential we get a finite graph $\Gamma(\rho) \subset \Sigma$ which we formally distinguish from the base set
$B_{\rho}$ (however as subsets in $\Sigma$ they are the same).  This graph carries  additional equipments ---  fixed points on each edge and the corresponding orinetation
for each segment ended at pairs of the closest marked points. 

The topological structure of the graph $\Gamma(\rho)$  is given by the following observation: by the very construction $\Gamma_{\rho}$ is homotopic to
$\Sigma \backslash \{p_1, ..., p_m \}$. Indeed, the base set is the union of critical points and finite trajectories, therefore
the punctured surface consists of the base set plus infinite trajectories, and we can shrink every infinite trajectory to point. We can compute 
the Euler characteristic $h_1(\Gamma(\rho)) - h_0(\Gamma(\rho))$ equals to $4 g - 4$ in the general case for holomorphic differential
without multiple zeros using the Mayer–Vietoris exact sequence, which is essensially
the same as the formula for the numbers of minima and saddle points above. For holomorphic differentials with multiple zeros one gets
the bound of the same type.  Thus the presence of $\Gamma(\rho)$ from the topological point of view doesn't distinguish different holomorphic differentials unless it has multiple zeros. 

But we can reduce some geometrical objects, f.e. vector bundles, from whole $\Sigma$ to the graphs which correspond to holomorphic differentials, and
al least we can define certain functions on the bundle ${\cal H} \to {\cal M}_g$ where ${\cal M}_g$ is the moduli space 
of Riemann surfaces and ${\cal H}$ is the  bundle whose fiber is $H^0(\Sigma, K_{\Sigma})$.  

As a toy example let us present a complex valued function on the total space $  {\cal H}$. To do this we first take holomorphic differentials
such that the corresponding functions $\Psi(\rho)$ have isolated finite critical points. In this situation each $\rho$ defines the graph
$\Gamma(\rho)$ and since the number of edges is finite one can correctly define the integral
$$
A_{SBS}: {\cal H} \to \mathbb{C}, \quad  A_{SBS}(\rho) = \int_{\Gamma (\rho)} \rho \in \mathbb{C},
$$
where the integration is done first on the oriented segments of $\Gamma(\rho)$ and then by summing up these numbers.

If $\Psi(\rho)$ admits degenerated critical points which form critical subsets then we still get the corresponding graph $\Gamma(\rho)$
admits certian edges with ``zero'' orientation being the components of the critical subset. Nevertheless the expression $A_{SBS}$ is
correctly defined in this situation: one just excludes the edges with zero orientation and integrates along the oriented segments only.

Then the natural probelm arises: what are the properties of the function $A_{SBS}$ on ${\cal H}$? Is it contineous? At least one thing
is clear from the defintion: it is linear along the fibers of ${\cal H} \to {\cal M}_g$. This means that it corresponds to a section
of the dual bundle ${\cal H}^* \to {\cal M}_g$. Take the zeroset of this section and denote it as $D_{SBS} \subset {\cal M}_g$.
What can be said about this subset? Is it a real submanifold? 

More serious questions arise when we see that all the story looks somehow related to a subject which is based on the same geometrical data: theory of flat surfaces, see
[4] and references therein. Every holomorphic differential $\rho$ defines a flat metric on $\Sigma$ with conical singularities at zeros
$P_1, ..., P_m$. One defines there certain dynamics on the total space ${\cal H}$ and on the components ${\cal H}(d_1, ..., d_m)$
where $d_i$ are multiplicities of zeros $P_1, ..., P_m$. It is natural to ask: what is the meaning of the function $A_{SBS}$ from
the point of view of flat surfaces?   

On the other hand the construction of finite graph $\Gamma (\rho)$ works for any positive degree of the canonical bundle: the Kahler structure
$(G, I, \Omega)$ induces the corresponding hermitian structure on $K_{\Sigma}^k$ for any $k \in \mathbb{Z}$, and the distinguished
prequantization connection $a_{LC}$ induces the corresponding connection $a_k \in {\cal A}_h(K^k_{\Sigma}$. The construction of $\Gamma(\rho_k)$
where $\rho_k \in H^0(\Sigma, K^k_{\Sigma})$ is holomorphic $k$ - differential is the same as in the case of $K_{\Sigma}$, therefore
we get certain new ingredients of any theory based on the geometry of the bundle ${\cal H}_k \to {\cal M}_g$, but it is already
non toy level.

$$$$

{\bf References:}

[1] {\bf N. Tyurin}, {\it ``Special Bohr - Sommerfeld geometry''}, arXiv: 1508.06804v1 [math.SG]; 

[2] {\bf P. Griffits, J. Harris}, {\it ``Principles of algebraic geometry''}, NY, Wiley, 1978;

[3] {\bf Y. Eliashberg}, {\it ``Topological characterization of Stein manifolds of $dim> 2$''}, Internat. J. Math., 1, no.1, pp. 29-46 (1990);

[4] {\bf A. Zorich}, {\it ``Flat surfaces''}, Frontiers in Number theory, Physics and geometry, Vol. 1, Springer 2006; arXiv: math/0609392v2.

\end{document}